\documentclass[11pt,reqno,a4paper,oneside]{amsart}

\usepackage{latexsym}
\usepackage{amsmath}
\usepackage{amsthm}
\usepackage{amssymb}
\usepackage{amsfonts}
\usepackage{comment}

\usepackage{mathrsfs}                           

\newcommand{\mR}{\mathbf{R}}                    
\newcommand{\mC}{\mathbf{C}}                    
\newcommand{\abs}[1]{\lvert #1 \rvert}          
\newcommand{\norm}[1]{\lVert #1 \rVert}         
\newcommand{\br}[1]{\langle #1 \rangle}         

\newcommand{\supp}{\mathrm{supp}}

\newcommand{\closure}[1]{\overline{#1}}

\newcommand{\mOp}{\mathrm{Op}}

\theoremstyle{definition}
\newtheorem{thm}{Theorem}[section]
\newtheorem{prop}[thm]{Proposition}

\newtheorem*{definition}{Definition}

\title[Inverse Problems for the Magnetic Schr\"odinger Equation]{Inverse Boundary Value Problems for the \\
Magnetic Schr\"odinger Equation}
\author{Mikko Salo}
\address{Department of Mathematics and Statistics, University of Helsinki, Finland}
\email{mikko.salo\@@helsinki.fi}
\date{}

\begin{document}

\begin{abstract}
We survey recent results on inverse boundary value problems for the magnetic Schr\"odinger equation.
\end{abstract}

\maketitle

\section{Introduction}

The purpose of this article is to give an overview of recent developments in inverse boundary value problems for the magnetic Schr\"odinger equation. In these problems, one attempts to find information about the coefficients of the equation from measurements made at the boundary. A prototype for this class of questions is the Calder\'on problem \cite{calderon}, also known as the inverse conductivity problem. There one wishes to determine the electrical conductivity of a body by measuring electrical currents corresponding to voltage potentials at the boundary. Both the theoretical and applied aspects of this problem have been extensively studied, see the surveys \cite{uhlmanndevelopments}, \cite{uhlmannselecta}.

Inverse problems for Schr\"o\-din\-ger equations are closely related to the Calder\'on problem. In fact, most results for the Calder\'on problem in three and higher dimensions are based on a reduction to the Schr\"odinger equation 
\begin{equation*}
(-\Delta + V)u  = 0.
\end{equation*}
Many other inverse problems, both for scalar equations and systems, have been treated using similar reductions. However, often these reductions will result in first order perturbations of the Laplacian, which correspond to equations of the form 
\begin{equation*}
(-\Delta + \sum_{j=1}^n A_j \frac{\partial}{\partial_{x_j}} + \tilde{V})u = 0.
\end{equation*}
The magnetic Schr\"odinger equation is a convenient model for these equations. By studying such first order perturbations, results have been obtained in inverse problems for isotropic elasticity \cite{nakamurauhlmann}, \cite{nakamurauhlmannerratum}, \cite{eskinralston_elasticity}, Maxwell equations in chiral media \cite{mcdowall}, Dirac equations \cite{nakamuratsuchida}, \cite{lidirac}, Schr\"odinger equations with external Yang-Mills potentials \cite{eskin}, and the Stokes system \cite{heckliwang}.

The magnetic Schr\"odinger equation is of course important in its own right, and it is a natural elliptic equation for which one can study inverse problems. It is well known that the inverse boundary value problem for the magnetic Schr\"odinger equation is equivalent to an inverse scattering problem at fixed energy, provided that the coefficients are compactly supported. Therefore, methods introduced for the inverse boundary value problem are often also useful in studying inverse scattering problems.

We wish to give a precise formulation of the inverse boundary problem. Let $\Omega \subseteq \mR^n$ be a bounded open set with Lipschitz boundary. We always assume $n \geq 3$ in this article. The magnetic Schr\"odinger operator is defined as 
\begin{equation*}
H_{A,V} = (D + A)^2 + V
\end{equation*}
where $A = (A_1,\ldots,A_n)$ is a vector field (magnetic potential) with components in $L^{\infty}(\Omega)$, and $V$ is a function (electric potential) in $L^{\infty}(\Omega)$. The coefficients are allowed to be complex valued. Also, $D = \frac{1}{i} \nabla$, and $\nabla + iA = i(D + A)$ is the magnetic gradient. In nondivergence form, $H_{A,V}$ is a first order perturbation of the Laplacian, given by 
\begin{equation*}
H_{A,V} = -\Delta + 2A \cdot D + \tilde{V},
\end{equation*}
where $\tilde{V} = A^2 + D \cdot A + V$.

Consider the weak Dirichlet problem 
\begin{equation*}
\left\{ \begin{array}{rll}
H_{A,V} u &\!\!\!= 0 & \quad \text{in } \Omega, \\
u &\!\!\!= f & \quad \text{on } \partial \Omega.
\end{array} \right.
\end{equation*}
By standard methods, $H_{A,V}$ has a countable set of Dirichlet eigenvalues in $\Omega$. We will make the standing assumption that $0$ is not a Dirichlet eigenvalue. In this case the Dirichlet problem has a unique weak solution $u = u_f \in H^1(\Omega)$ for any $f \in H^{1/2}(\partial \Omega)$.

The boundary measurements related to $H_{A,V}$ are given by the Dirichlet-to-Neumann map (DN map) 
\begin{equation*}
\Lambda_{A,V}: f \mapsto (\nabla + iA)u_f \cdot \nu|_{\partial \Omega}.
\end{equation*}
Here $\nu$ is the outer unit normal of $\partial \Omega$. Thus, for any boundary value $f$, one can measure the magnetic normal derivative of $u_f$ at the boundary. More precisely, we define $\Lambda_{A,V}$ using the weak formulation 
\begin{equation} \label{dnmap_def}
(\Lambda_{A,V} f | g)_{\partial \Omega} = ((\nabla+iA) u_f | (\nabla+i\bar{A}) e_g) + (Vu_f|e_g), \quad f, g \in H^{1/2}(\partial \Omega),
\end{equation}
where $e_g$ is any function in $H^1(\Omega)$ with $e_g|_{\partial \Omega} = g$, and where we write 
\begin{equation*}
(u|v) = \int_{\Omega} u \cdot \bar{v} \,dx, \quad (f|g)_{\partial \Omega} = \int_{\partial \Omega} f \bar{g} \,dS.
\end{equation*}
Here $dS$ is the Euclidean surface measure on $\partial \Omega$. It follows from \eqref{dnmap_def} that $\Lambda_{A,V}$ is a bounded map $H^{1/2}(\partial \Omega) \to H^{-1/2}(\partial \Omega)$.

An important property of the magnetic Schr\"odinger operator is gauge invariance. This means that in the weak sense 
\begin{equation} \label{h_gauge}
e^{-i\varphi} H_{A,V}(e^{i\varphi} u) = H_{A+\nabla \varphi,V} u
\end{equation}
if $\varphi \in W^{1,\infty}(\Omega)$ and $u \in H^1(\Omega)$. If additionally $\varphi|_{\partial \Omega} = 0$, the boundary measurements are preserved:
\begin{equation} \label{dn_gauge}
\Lambda_{A+\nabla \varphi,V} = \Lambda_{A,V}.
\end{equation}
This gauge transformation replaces the magnetic potential $A$ by a gauge equivalent potential $A + \nabla \varphi$. However, the magnetic field $dA$ is preserved, where $dA = d(\sum_{j=1}^n A_j(x) \,dx_j)$ is the exterior derivative of $A$ interpreted as a $1$-form. If $n = 3$, then $dA$ corresponds to $\nabla \times A$.

We may now formulate the basic inverse problem considered in this article.
\begin{quote}
{\bf Inverse problem for the magnetic Schr\"odinger equation}: Given the boundary measurements $\Lambda_{A,V}$, is it possible to determine the magnetic field $dA$ and the electric potential $V$ in $\Omega$?.
\end{quote}
As with other inverse problems, there are different aspects which can be studied. We will concentrate on the following topics.
\begin{enumerate}
\item 
{\bf Boundary determination}: given $\Lambda_{A,V}$, determine the values of $A$ and $V$ on $\partial \Omega$.
\item 
{\bf Uniqueness}: if $\Lambda_{A_1,V_1} = \Lambda_{A_2,V_2}$, show that $dA_1 = dA_2$ and $V_1 = V_2$.
\item 
{\bf Reconstruction}: given $\Lambda_{A,V}$, reconstruct the values of $dA$ and $V$.
\item 
{\bf Stability}: if $\Lambda_{A_1,V_1}$ and $\Lambda_{A_2,V_2}$ are close, show that $dA_1$ and $dA_2$ (and $V_1$ and $V_2$) also have to be close.
\item 
{\bf Partial data}: if $\Lambda_{A_1,V_1}|_{\Gamma} = \Lambda_{A_2,V_2}|_{\Gamma}$ for a subset $\Gamma \subseteq \partial \Omega$, show that $dA_1 = dA_2$ and $V_1 = V_2$.
\end{enumerate}

The inverse problem for the magnetic Schr\"odinger equation has been studied by several authors. One of the first results is due to Sun \cite{sun}, who showed uniqueness for magnetic potentials in $W^{2,\infty}$ satisfying a smallness condition. The proof was based on complex geometrical optics (CGO) solutions, which were introduced by Sylvester and Uhlmann \cite{sylvesteruhlmann} in the case $A = 0$. If $A$ is nonzero and large, then the construction of CGO solutions is more difficult, and was achieved by Nakamura and Uhlmann \cite{nakamurauhlmann}, \cite{nakamurauhlmannerratum} using a pseudodifferential conjugation method. With this method, Nakamura, Sun, and Uhlmann \cite{nakamurasunuhlmann} showed uniqueness for smooth coefficients and domains without any smallness conditions, and they also gave a boundary determination result. The uniqueness result was extended to $C^1$ magnetic potentials by Tolmasky \cite{tolmasky}, and to some less regular but small potentials by Panchenko \cite{panchenko}.

The fixed energy inverse scattering problem for magnetic Schr\"odinger operators has been studied by Henkin and Novikov \cite{novikovkhenkin} for small $A$, and by Eskin and Ralston \cite{eskinralston} without smallness assumptions. There is a close relation between the inverse boundary value problem and the fixed energy inverse scattering problem, and in fact for compactly supported potentials the two problems are equivalent. See Eskin and Ralston \cite{eskinralstonproc} for a proof of this in the magnetic case. We will not consider scattering problems in this article.

Recently, there have been many developments in the above questions. Brown and Salo \cite{brownsalo} proved a boundary determination result for continuous magnetic potentials in $C^1$ domains. Interior uniqueness was shown to hold for Dini continuous magnetic potentials by Salo \cite{salothesis}. Salo \cite{saloreconstruction} gave a reconstruction algorithm, and stability results were proved in Tzou \cite{tzou}. A partial data result was given by Dos Santos Ferreira, Kenig, Sj\"ostrand, and Uhlmann \cite{dksu} for $C^2$ magnetic potentials, following the result for $A = 0$ by Kenig, Sj\"ostrand, and Uhlmann \cite{ksu}. These works also introduced a new method for constructing CGO solutions, based on Carleman estimates and a convexification idea. The partial data result was extended to H\"older continuous potentials by Knudsen and Salo \cite{knudsensalo}.

We will describe these recent developments below. The article is organized as follows. Section \ref{sec:boundary} explains the boundary determination result in \cite{brownsalo}. In Section \ref{sec:cgo_pdo}, we discuss the construction of CGO solutions using pseudodifferential conjugation, following \cite{saloreconstruction}. Section \ref{sec:uniqueness} is devoted to the uniqueness, reconstruction, and stability results given in \cite{salothesis}, \cite{saloreconstruction}, \cite{tzou}. In Section \ref{sec:cgo_carleman} we describe the construction of CGO solutions based on Carleman estimates following \cite{dksu} and \cite{knudsensalo}, and the final Section \ref{sec:partial} considers the partial data results given in these papers.

\vspace{5pt}

\noindent {\bf Notation.} All functions will be complex valued, unless stated otherwise. We let $C^k(\closure{\Omega})$ be the set of functions whose partial derivatives up to order $k$ are continuous in $\closure{\Omega}$, and $C^{k+\varepsilon}(\closure{\Omega})$ consists of the functions in $C^k(\closure{\Omega})$ whose partial derivatives of order $k$ are $\varepsilon$-H\"older continuous in $\closure{\Omega}$. The set of compactly supported continuous functions in $\mR^n$ is denoted by $C_c(\mR^n)$, similarly $C^{\infty}_c(\mR^n)$ and $L^p_c(\mR^n)$. The $L^2$-based and $L^p$-based Sobolev spaces are denoted by $H^s$ and $W^{k,p}$, respectively, and $H^1_0(\Omega)$ is the set of functions in $H^1(\Omega)$ with zero boundary values. Functions with values in a Banach space $X$ are denoted by $C^k(\closure{\Omega} ; X)$ etc.

\section{Boundary determination} \label{sec:boundary}

In the context of the Calder\'on problem, if the domain and conductivity are smooth, it was shown in \cite{sylvesteruhlmannboundary} that the DN map is a pseudodifferential operator, and the Taylor series of the conductivity at the boundary can be read off from the symbol. The same philosophy applies to many inverse problems. For the magnetic Schr\"odinger operator, this method was used in \cite{nakamurasunuhlmann} to show that if everything is smooth, then the Taylor series of the tangential components of $A$ and of $V$ at the boundary can be determined from $\Lambda_{A,V}$ (there is a small mistake in \cite[Theorem D]{nakamurasunuhlmann}, but this is not hard to fix). The tangential components of $A$ on $\partial \Omega$ are 
\begin{equation} \label{tangential_component}
A_{\mathrm{tan}} = A - (A \cdot \nu)\nu.
\end{equation}
The gauge invariance \eqref{dn_gauge} shows that one can only expect to recover $A_{\mathrm{tan}}$ from $\Lambda_{A,V}$.

If the domain and coefficients are not smooth, then pseudodifferential theory is not available. However, the symbol of a pseudodifferential operator can be obtained by testing against oscillatory functions, and this idea works also in the nonsmooth case. Early results are given in \cite{kohnvogelius1}. For the Calder\'on problem in Lipschitz domains the method was used in \cite{brownboundary}, and the magnetic Schr\"odinger equation in $C^2$ domains was  considered in \cite{salothesis}. The following theorem is proved in \cite{brownsalo}.

\begin{thm} \label{thm:boundary}
Let $\Omega$ be a $C^1$ domain and $n \geq 3$. Let $A \in C(\closure{\Omega} ; \mC^n)$ and $V \in L^{\infty}(\Omega)$. Then $\Lambda_{A,V}$ determines $A_{\mathrm{tan}}$ on $\partial \Omega$.
\end{thm}

More precisely, let $x_0 \in \partial \Omega$ and let $\alpha$ be a unit tangent vector to $\partial \Omega$ at $x_0$. In \cite{brownsalo} one constructs a family of oscillating functions $f_M \in L^2(\partial \Omega)$, with $f_M$ supported in $B(x_0,1/M) \cap \partial \Omega$, such that 
\begin{equation} \label{boundary_limit}
\lim_{M \to \infty} ((\Lambda_{A,V}-\Lambda_{0,0})f_M|f_M)_{\partial \Omega} = \alpha \cdot A(x_0).
\end{equation}
The method is local: one only needs to know $\Lambda_{A,V} f$ near $x_0$ for functions $f$ supported near $x_0$. The method also gives the stability estimate 
\begin{equation*}
\norm{(A_1-A_2)_{\mathrm{tan}}}_{L^{\infty}(\partial \Omega)} \leq \sqrt{2} \norm{\Lambda_{A_1,V_1} - \Lambda_{A_2,V_2}}_{L^2(\partial \Omega) \to L^2(\partial \Omega)},
\end{equation*}
if $A_j$ are continuous in $\closure{\Omega}$ and $V_j$ are $L^{\infty}$.

We now describe some ideas in the proof. Let $\Omega$ be given by a defining function $\rho \in C^1(\mR^n ; \mR)$, so that $\Omega = \{ \rho > 0 \}$ and $\partial \Omega = \{ \rho = 0 \}$, and $\nabla \rho \neq 0$ on $\partial \Omega$. Let $x_0$ and $\alpha$ be as above, and normalize $\rho$ so that $\nabla \rho(x_0) = -\nu(x_0)$. We also choose coordinates so that $x_0 = 0$ and $\nu(x_0) = -e_n$, the $n$th coordinate vector.

Suppose $\omega$ is a modulus of continuity for $\nabla \rho$, and let $\eta \in C^{\infty}_c(\mR^n)$ be a function supported in $B(0,1/2)$ and satisfying $\int_{\mR^{n-1}} \eta(x',0)^2 \,dx' = 1$. Put $\eta_M(x) = \eta(M(x',\rho(x)))$ where $M$ is a large parameter, and define 
\begin{eqnarray*}
 & u_0(x) = \eta_M(x) e^{N(i\alpha \cdot x - \rho(x))}, & x \in \mR^n, \\
 & f_M(x) = M^{\frac{n-1}{2}} u_0(x), & x \in \partial \Omega,
\end{eqnarray*}
where $N$ is a parameter larger than $M$, defined by 
\begin{equation*}
M^{-1} \omega(M^{-1}) = N^{-1}.
\end{equation*}
The idea is that $u_0$ is an approximate solution to $H_{A,V} u = 0$ in $\Omega$, concentrated near the boundary point $0$ and oscillating at the boundary.

From \eqref{dnmap_def} we obtain 
\begin{multline} \label{boundary_mainid}
((\Lambda_{A,V}-\Lambda_{0,0})f_M|f_M)_{\partial \Omega} = M^{n-1} \int_{\Omega} (iA \cdot (u\nabla\bar{v} - \bar{v}\nabla u) + (A^2 + V) u\bar{v}) \,dx
\end{multline}
where $u \in H^1(\Omega)$ solves $H_{A,V} u = 0$ with $u = u_0$ on $\partial \Omega$, and $v \in H^1(\Omega)$ solves $-\Delta v = 0$ with $v = u_0$ on $\partial \Omega$. We write $u = u_0 + u_1$ and $v = u_0 + v_1$ where $u_1, v_1 \in H^1_0(\Omega)$. Since $u_0$ was chosen as an approximate solution of $H_{A,V} u = 0$ (and $-\Delta u = 0$), $u_1$ and $v_1$ will be error terms smaller than $u_0$, which solve 
\begin{eqnarray}
 & H_{A,V} u_1 = -H_{A,V} u_0 & \quad \text{in } \Omega, \label{u1_eq} \\
 & -\Delta v_1 = \Delta u_0 & \quad \text{in } \Omega. \label{v1_eq}
\end{eqnarray}

Using the explicit form of $u_0$, the main term in \eqref{boundary_mainid} will satisfy 
\begin{equation*}
\lim_{M \to \infty} M^{n-1} \int_{\Omega} iA \cdot (u_0\nabla\bar{u}_0 - \bar{u}_0 \nabla u_0 \,dx = A(0) \cdot \alpha.
\end{equation*}
Thus, \eqref{boundary_limit} and Theorem \ref{thm:boundary} will follow if one can show that all the terms in \eqref{boundary_mainid} involving $u_1$ or $v_1$ will be $o(1)$ as $M \to \infty$.

Since $u_1, v_1$ are $H^1_0(\Omega)$ solutions of \eqref{u1_eq}, \eqref{v1_eq} where the right hand sides are explicit, we get good bounds for $\norm{\nabla u_1}_{L^2(\Omega)}$ and $\norm{\nabla v_1}_{L^2(\Omega)}$ from weak solution estimates. The problem is to estimate $\norm{u_1}_{L^2(\Omega)}$ and $\norm{v_1}_{L^2(\Omega)}$. We will use Hardy's inequality, 
\begin{equation*}
\norm{w/\delta}_{L^2(\Omega)} \leq C \norm{\nabla w}_{L^2(\Omega)}, \qquad w \in H^1_0(\Omega),
\end{equation*}
where $\delta(x) = \mathrm{dist}(x,\partial \Omega)$, and also the estimate 
\begin{equation*}
\norm{\delta \nabla w}_{L^2(\Omega)} \leq C \norm{w}_{L^2(\Omega)}, \qquad w \in L^2(\Omega) \text{ harmonic},
\end{equation*}
which follows from $L^2$ boundedness of the Hardy-Littlewood maximal function (see \cite[Lemma 5.5]{salothesis}). Using these two estimates in \eqref{boundary_mainid}, it will be enough to obtain a bound for $\norm{v_1}_{L^2(\Omega)}$.

The main technical point of \cite{brownsalo} is to prove for $C^1$ domains the estimate 
\begin{equation} \label{v1_estimate}
\norm{v_1}_{L^2(\Omega)} \leq C M^{\frac{1-n}{2}} N^{-\frac{1}{2}}.
\end{equation}
This will allow to finish the proof of Theorem \ref{thm:boundary}. To prove \eqref{v1_estimate} we first note that the Poincar\'e inequality and the gradient bound for $v_1$ yield a good estimate for $\norm{v_1}_{L^2(\Omega \cap B(0,100M^{-1}))}$. In $\Omega \smallsetminus B(0,100M^{-1})$ we know that $v_1$ is a harmonic function vanishing on $\partial \Omega$. Now we use the fact that there is another such function for which good bounds are known, namely the Green function $G(x,y)$ for the Laplacian in $\Omega$, with singularity at a suitable point in $B(0,100M^{-1})$. The pointwise values of $v_1$ can be compared with $G$, using the boundary Harnack principle of Jerison and Kenig \cite{jkexp}. In high dimensions, one gets a bound for $\norm{v_1}_{L^2(\Omega \smallsetminus B(0,100M^{-1}))}$ from the easy estimate 
\begin{equation*}
\abs{G(x,y)} \leq C \abs{x-y}^{2-n}.
\end{equation*}
In lower dimensions one needs to use additional vanishing of $G$ at the boundary. We refer to \cite{brownsalo} for the details.

\section{CGO solutions and pseudodifferential conjugation} \label{sec:cgo_pdo}

In \cite{sylvesteruhlmann}, CGO solutions were introduced for the equation $H_{0,V} u = 0$. They are of the form 
\begin{equation} \label{cgo_def_1}
u = e^{i\rho \cdot x}(1+r_{\rho})
\end{equation}
where $\rho \in \mC^n$ is a complex vector satisfying $\rho \cdot \rho = 0$, and $r_{\rho}$ is a remainder term which is small when $\abs{\rho}$ is large. Thus, for $\abs{\rho}$ large the solution $u$ resembles the harmonic exponential $e^{i\rho \cdot x}$, and one can use these exponentials to recover the potential $V$ from boundary measurements $\Lambda_{0,V}$. The method of CGO solutions is central in most results on inverse boundary problems for elliptic equations.

We wish to discuss CGO solutions to the magnetic Schr\"odinger equation $H_{A,V} u = 0$, or equivalently to $(-\Delta + 2A \cdot D + \tilde{V}) u = 0$. This is more difficult than for $H_{0,V}$. The first construction was achieved by Nakamura and Uhlmann \cite{nakamurauhlmann}, \cite{nakamurauhlmannerratum}, who introduced a pseudodifferential conjugation argument which essentially removes $A$. Stronger forms of this construction are now available. In \cite{tolmasky} the argument was extended to $C^{2/3+\varepsilon}$ magnetic potentials, and in \cite{salothesis} to just continuous $A$. These results only yield solutions in a bounded domain, and there is no uniqueness for the solutions. In \cite{saloreconstruction} the construction was extended to yield global unique CGO solutions.

We will follow \cite{saloreconstruction}. First we introduce weighted $L^2$ spaces $L^2_{\delta}(\mR^n)$ with norm $\norm{\br{x}^{\delta} f}_{L^2(\mR^n)}$, where $\br{x} = (1+\abs{x}^2)^{1/2}$, and weighted Sobolev spaces $H^k_{\delta}(\mR^n)$ with norm $\norm{f}_{H^k_{\delta}} = \sum_{\abs{\alpha} \leq k} \norm{\partial^{\alpha} f}_{L^2_{\delta}}$. The following is the main result in this section. It shows existence and uniqueness of CGO solutions for large $\abs{\rho}$ with known asymptotics as $\abs{\rho} \to \infty$.

\begin{thm} \label{thm:cgo_pdo}
Suppose $A \in C_c(\mR^n ; \mC^n)$ and $\tilde{V} \in L^{\infty}_c(\mR^n)$. Let $-1<\delta<0$. If $\rho \in \mC^n$ satisfies $\rho \cdot \rho = 0$ and if $\abs{\rho}$ is large enough, then the equation 
\begin{equation*}
(-\Delta + 2A \cdot D + \tilde{V}) u = 0 \quad \text{in } \mR^n
\end{equation*}
has a unique CGO solution of the form \eqref{cgo_def_1} with $r_{\rho} \in H^1_{\delta}$. Further, if $\rho$ has the form $\rho = \frac{1}{h}(\nu_1+i\nu_2)$ where $\nu_1,\nu_2 \in \mR^n$ are orthogonal unit vectors, and 
\begin{equation} \label{asymptotic_phi_def}
\phi(x) = -\frac{1}{2\pi} \int_{\mR^2} \frac{1}{y_1+i y_2} (\nu_1+i\nu_2) \cdot A(x-y_1 \nu_1-y_2 \nu_2) \,dy_1 \,dy_2,
\end{equation}
then one has the asymptotics $1 + r_{\rho} = a_{\rho} + \tilde{r}_{\rho}$ where 
\begin{eqnarray}
 & a_{\rho} \to e^{i \phi} \text{ pointwise}, \quad \norm{a_{\rho}}_{L^{\infty}} = O(1), \quad h \norm{\nabla a_{\rho}}_{L^{\infty}} = o(1), & \label{arho_asymptotics} \\
 & \norm{\tilde{r}_{\rho}}_{L^2_{\delta}} + h \norm{\nabla \tilde{r}_{\rho}}_{L^2_{\delta}} = o(1), & \label{rtilderho_asymptotics}
\end{eqnarray}
as $h \to 0$.
\end{thm}

To discuss the proof of Theorem \ref{thm:cgo_pdo}, we introduce the conjugated operator 
\begin{equation*}
e^{-i\rho \cdot x} H_{A,V} e^{i\rho \cdot x} = P_{\rho} + 2A \cdot D_{\rho} + \tilde{V}
\end{equation*}
where $P_{\rho} = -\Delta + 2\rho \cdot D$, $D_{\rho} = D + \rho$, and $\tilde{V} = A^2 + D \cdot A + V$. Inserting \eqref{cgo_def_1} in $H_{A,V} u = 0$, we see that the construction of CGO solutions reduces to solving 
\begin{equation} \label{inhomog_eq}
(P_{\rho} + 2A \cdot D_{\rho} + \tilde{V})r_{\rho} = f
\end{equation}
where $f = - 2A \cdot \rho - \tilde{V}$. If $A = V = 0$ the solvability of \eqref{inhomog_eq} follows from the fundamental estimates of Sylvester and Uhlmann.

\begin{prop} \label{prop:sylvesteruhlmann}
\cite{sylvesteruhlmann} Let $-1 < \delta < 0$, and let $\rho \in \mC^n$ with $\rho \cdot \rho = \lambda$ and $\abs{\rho} \geq 1$. Then for any $f \in L^2_{\delta+1}$ the equation $P_{\rho} u = f$ has a unique solution $u \in L^2_{\delta}$. The solution operator $G_{\rho}: f \mapsto u$ satisfies 
\begin{equation*}
\abs{\rho} \norm{G_{\rho} f}_{L^2_{\delta}} + \norm{\nabla G_{\rho} f}_{L^2_{\delta}} \leq C \norm{f}_{L^2_{\delta+1}}.
\end{equation*}
\end{prop}

If $\norm{A}_{L^{\infty}}$ is small then we may solve \eqref{inhomog_eq} by trying $r_{\rho} = G_{\rho} v_{\rho}$ and using Neumann series. If $A$ is large this will not work, since then $I + 2A \cdot D_{\rho} G_{\rho}$ is not a small perturbation of identity on $L^2_{\delta+1}$.

If $A$ is large, one might try to use a gauge transformation as in \eqref{h_gauge} to remove $A$. However, if $A$ is not a gradient, this will not remove all of $A$. The idea of Nakamura and Uhlmann was to replace the exponentials $e^{i\varphi}$ in \eqref{h_gauge} by more general pseudodifferential operators. This may be thought of as a pseudodifferential gauge transformation, which essentially removes $A$. The rest of the proof then proceeds by a perturbation argument as in the case where $A$ is small, using the norm estimates for $G_{\rho}$ and Neumann series.

It will be convenient to switch to semiclassical notation. This amounts to writing $h = \frac{1}{\abs{\rho}}$ and to considering a small parameter $h$ instead of a large parameter $\abs{\rho}$. However, we may then use the well-established machinery of semiclassical pseudodifferential calculus, which automatically keeps track of the dependence on $h$ in norm estimates, and simplifies some proofs since the parameter $h$ is scaled away when passing to symbols.

We will use the usual semiclassical symbol classes, see \cite{dimassisjostrand}.
\begin{definition}
If $0 \leq \sigma < 1/2$ and $m \in \mR$, we let $S^m_{\sigma}$ be the space of all functions $c(x,\xi) = c(x,\xi;h)$ where $x,\xi \in \mR^n$ and $h \in (0,h_0]$, $h_0 \leq 1$, such that $c$ is smooth in $x$ and $\xi$ and 
\begin{equation*}
\abs{\partial_x^{\alpha} \partial_{\xi}^{\beta} c(x,\xi)} \leq C_{\alpha \beta} h^{-\sigma\abs{\alpha+\beta}} \br{\xi}^m
\end{equation*}
for all $\alpha,\beta$. If $c \in S^m_{\sigma}$ we define an operator $C = \mOp_h(c)$ by 
\begin{equation*}
Cf(x) = (2\pi)^{-n} \iint_{\mR^{2n}} e^{i(x-y)\cdot\xi} c\left(\frac{x+y}{2},h\xi\right) f(y) \,d\xi \,dy.
\end{equation*}
\end{definition}

Semiclassical quantization means that the symbol $\xi_j$ corresponds to the operator $hD_j$. Note that we have used Weyl quantization. The operators have a calculus: if $c \in S^0_{\sigma}$ then $C$ is bounded on weighted spaces $L^2_{\delta}$ with norm independent of $h$, and if $c \in S^m_{\sigma}$ and $d \in S^{m'}_{\sigma}$ then $CD$ is an operator with symbol in $S^{m+m'}_{\sigma}$.

To manage the nonsmooth coefficients, we introduce the standard mollifier $\chi_{\delta}(x) = \delta^{-n} \chi(x/\delta)$ where $\chi \in C_c^{\infty}(\mR^n)$, $0 \leq \chi \leq 1$, and $\int \chi \,dx = 1$. Fix $\sigma_0$ with $0 < \sigma_0 < 1/2$, and consider the $h$-dependent smooth approximations 
\begin{equation} \label{asharp_def}
A^{\sharp} = A \ast \chi_{\delta},
\end{equation}
with the specific choice $\delta = h^{\sigma_0}$.

We may write the equation \eqref{inhomog_eq} in semiclassical notation as 
\begin{equation}  \label{conjugated_identity_pdo}
(P + hQ + h^2 \tilde{V}) r_{\rho} = h^2 f
\end{equation}
where $P = -h^2 \Delta + 2 \hat{\rho} \cdot hD$ and $Q = 2A \cdot (hD + \hat{\rho})$, and $\hat{\rho} = h\rho$ so $\abs{\hat{\rho}} = 1$. We split $Q$ as $Q = Q^{\sharp} + Q^{\flat}$ where $Q^{\sharp} = 2A^{\sharp} \cdot (hD + \hat{\rho})$. The following result states that the smooth first order term $Q^{\sharp}$ can be conjugated into a term which has order less than $1$.

\begin{prop} \label{prop:conjugation}
Given $\sigma$ with $\sigma_0 < \sigma < 1/2$, there exist $c,\tilde{c},s \in S^0_{\sigma}$ and $\delta > 0$ such that 
\begin{equation} \label{conjugationidentity}
(P + hQ^{\sharp})C = \tilde{C}P + h^{1+\delta} \br{x}^{-1} S.
\end{equation}
Further, $C$ and $\tilde{C}$ are elliptic, in the sense that $c$ and $\tilde{c}$ are nonvanishing for small $h$.
\end{prop}

To prove Proposition \ref{prop:conjugation}, one takes some $c \in S^0_{\sigma}$ and commutes $P + hQ^{\sharp}$ and $C$ using the symbol calculus. The result is 
\begin{equation*}
(P + hQ^{\sharp})C = CP + h \mOp_h(\frac{1}{i} H_p c + q^{\sharp}c) + \text{lower order terms}.
\end{equation*}
Here $p(\xi) = \xi^2 + 2\hat{\rho} \cdot \xi$ is the symbol of $P$, $q^{\sharp} = 2A^{\sharp} \cdot (\xi + i\nabla\varphi) - hD \cdot A^{\sharp}$, and $H_p = \nabla_{\xi} p \cdot \nabla_x - \nabla_x p \cdot \nabla_{\xi} = 2(\xi+\hat{\rho}) \cdot \nabla_x$ is the Hamilton vector field. We want to choose $c$ so that the first order term disappears, i.e. 
\begin{equation} \label{hp_dbar_eq}
\frac{1}{i} H_p c + q^{\sharp} c \approx 0.
\end{equation}
Trying $c = e^{i\varphi}$, this equation becomes 
\begin{equation*}
H_p \varphi \approx -q^{\sharp}.
\end{equation*}
Near $p^{-1}(0)$ this looks like a $\bar{\partial}$-equation in some variables, but degenerates away from $p^{-1}(0)$. However, the operator $P$ is elliptic away from $p^{-1}(0)$, and it will be sufficient to find $\varphi$ near $p^{-1}(0)$. Inserting a suitable cutoff and applying the Cauchy transform, we obtain a function $\varphi$ satisfying 
\begin{equation*}
\abs{\partial_x^{\alpha} \partial_{\xi}^{\beta} \varphi(x,\xi)} \leq C_{\alpha \beta} h^{-\sigma_0 \abs{\alpha+\beta}} \br{x}^{\abs{\beta}-1}.
\end{equation*}
The problem is that the $\xi$-derivatives of $\varphi$ grow in $x$, so $\varphi$ is not in $S^0_{\sigma_0}$. This can be fixed by introducing another cutoff $\chi \in C^{\infty}_c(\mR^n)$ with $\chi = 1$ near $0$, and by taking 
\begin{equation*}
c = e^{i\chi(h^{\theta} x)\varphi}
\end{equation*}
where $\theta > 0$ is chosen so that $\sigma := \sigma_0 + \theta < 1/2$. With this choice $c \in S^0_{\sigma}$, and we have traded off the undesired $x$-growth into slightly worse $h$-dependence in the symbol. This method was used in the context of energy estimates for nonlinear Schr\"odinger equations by Kenig, Ponce, and Vega \cite{kenigponcevega}.

We have arrived at the correct conjugating operator $C$, and it is possible to finish the proof of Proposition \ref{prop:conjugation}. Theorem \ref{thm:cgo_pdo} follows by combining the pseudodifferential conjugation with a perturbation argument where small errors are inverted by Neumann series. The asymptotics for CGO solutions may be computed either from the asymptotics of $C$, or by guessing the right asymptotics and using the uniqueness part of Theorem \ref{thm:cgo_pdo}. The details are given in \cite{saloreconstruction}.

\section{Uniqueness, reconstruction, and stability} \label{sec:uniqueness}

The CGO solutions constructed in Theorem \ref{thm:cgo_pdo} can be used to solve the uniqueness, reconstruction, and stability questions for the inverse problem for $H_{A,V}$. We begin with reconstruction. The next result was proved in \cite{saloreconstruction}. We state it here in a simplified form to avoid technicalities.

\begin{thm} \label{thm:reconstruction}
Let $\Omega \subseteq \mR^n$, $n \geq 3$, be a bounded simply connected domain with smooth boundary. Suppose $A \in C^{1+\varepsilon}_c(\Omega ; \mC^n)$ and $V \in L^{\infty}(\Omega ; \mC)$. Then one may reconstruct $dA$ and $V$ from $\Lambda_{A,V}$.
\end{thm}

In the case $A = 0$ this result is due to Nachman \cite{nachman}, who used CGO solutions and ideas from scattering theory to reconstruct $V$ from $\Lambda_{0,V}$. The fact that one has unique global CGO solutions for $H_{0,V}$ is crucial for the algorithm. Theorem \ref{thm:cgo_pdo} provides such solutions for $H_{A,V}$, and using these we can extend Nachman's algorithm to the magnetic case.

If $\rho \in \mC^n$ satisfies $\rho \cdot \rho = 0$ and $\abs{\rho}$ is large, we denote by $u_{\rho}$ the unique CGO solution to $H_{A,V} u = 0$ given by Theorem \ref{thm:cgo_pdo}. The first step is to show that $u_{\rho}$ may be characterized in terms of an integral equation on $\partial \Omega$. More precisely, $u_{\rho}|_{\partial \Omega}$ is the unique solution in $H^{3/2}(\partial \Omega)$ to the equation 
\begin{equation} \label{boundary_integral_eq}
(I + S_{\rho}(\Lambda_{A,V} - \Lambda_{0,0}))f = e^{i\rho \cdot x} \quad \text{on } \partial \Omega,
\end{equation}
where $S_{\rho}$ is the single layer potential 
\begin{equation*}
S_{\rho} f(x) = \int_{\partial \Omega} G_{\rho}(x,y) f(y) \,dS(y), \quad f \in H^{1/2}(\partial \Omega), x \in \partial \Omega,
\end{equation*}
and $G_{\rho}$ is the $\rho$-dependent Green function 
\begin{equation*}
G_{\rho}(x,y) = e^{i\rho \cdot (x-y)} (2\pi)^{-n} \int \frac{e^{i(x-y) \cdot \xi}}{\xi^2 + 2\rho \cdot \xi} \,d\xi.
\end{equation*}
Note that if $\rho = 0$ this would be the usual Green function for $-\Delta$.

The point is that if we know the boundary measurements $\Lambda_{A,V}$, and if $\rho$ is fixed, then we know the quantities in \eqref{boundary_integral_eq} and can determine $u_{\rho}|_{\partial \Omega}$ as the unique solution of this equation. Then we also know the scattering transform $t_{A,V}$, defined by 
\begin{equation*}
t_{A,V}(\xi,\rho) = ((\Lambda_{A,V}-\Lambda_{0,-\abs{\xi}^2})(u_{\rho}|_{\partial \Omega}), e^{ix \cdot (\xi+\rho)})_{\partial \Omega},
\end{equation*}
where $\xi \in \mR^n$ is a vector, and $\rho \in \mC^n$ satisfies $\rho \cdot \rho = 0$, $\mathrm{Re}\,\rho \perp \xi$, $\mathrm{Im}\,\rho \perp \xi$, and $\abs{\rho}$ is large.

The expression $t_{A,V}$ is chosen so that in the limit $\abs{\rho} \to \infty$, we obtain certain nonlinear Fourier transforms of components of $dA$. This motivates the name ''scattering transform'', since similar phenomena occur in scattering theory. In fact, using the definition \eqref{dnmap_def} of $\Lambda_{A,V}$ and the asymptotics for $u_{\rho}$ given in Theorem \ref{thm:cgo_pdo}, one sees that 
\begin{equation} \label{htav_limit}
\lim_{h \to 0} h\,t_{A,V}(\xi,\frac{1}{h}(\nu_1+i\nu_2)) = 2 \int e^{-ix\cdot\xi} e^{i\phi} (\nu_1+i\nu_2) \cdot A \,dx
\end{equation}
where $\xi,\nu_1,\nu_2$ are orthogonal and $\abs{\nu_1} = \abs{\nu_2} = 1$, and $\phi$ is as in \eqref{asymptotic_phi_def}. If one replaces $e^{i\phi}$ by $1$ this would just be the Fourier transform of a component of $dA$, so this indeed looks like a nonlinear Fourier transform of $dA$.

However, due to a result in \cite{eskinralston}, this purported nonlinear Fourier transform is just the usual one: integration by parts shows 
\begin{equation*}
\int e^{-ix\cdot\xi} e^{i\phi} (\nu_1+i\nu_2) \cdot A \,dx = \int e^{-ix\cdot\xi} (\nu_1+i\nu_2) \cdot A \,dx.
\end{equation*}
Thus, given $\Lambda_{A,V}$, we have recovered the last expression for different orthogonal triplets $\xi,\nu_1,\nu_2$ where $\abs{\nu_1} = \abs{\nu_2} = 1$. These are just the Fourier transforms of components of $dA$, and we have reconstructed the magnetic field.

After finding $dA$ we want to reconstruct $V$. First we construct some $C^{1+\varepsilon}$ magnetic potential $\tilde{A}$ with $d\tilde{A} = dA$ and $\tilde{A}|_{\partial \Omega} = 0$. Since $\Omega$ is simply connected, we know that $\tilde{A} = A + \nabla p$ where $p|_{\partial \Omega} = 0$, possibly after substracting a constant. By gauge equivalence, we know $\Lambda_{A,0} = \Lambda_{\tilde{A},0}$. For $\xi \in \mR^n$, define another scattering transform 
\begin{equation*}
\tilde{t}_{A,V}(\xi) = ((\Lambda_{A,V}-\Lambda_{A,0})(u_{\rho}|_{\partial \Omega}), v_{\tilde{\rho}}|_{\partial \Omega})_{\partial \Omega}
\end{equation*}
where $\rho$, $\tilde{\rho}$ are complex vectors chosen in a certain way so that $\rho + \tilde{\rho} = -\xi$ and $\abs{\rho} = \abs{\tilde{\rho}} = 1/h$, and $u_{\rho}$ and $v_{\tilde{\rho}}$ are suitable CGO solutions whose boundary values are obtained as solutions of integral equations. It follows that 
\begin{equation*}
\lim_{h \to \infty} \tilde{t}_{A,V}(\xi) = \int e^{-ix\cdot\xi} V(x) \,dx,
\end{equation*}
and the electric potential is determined.

If one is only interested in uniqueness, the regularity conditions can be relaxed somewhat. We say that $f$ is Dini continuous, written $f \in C^{\mathrm{Dini}}(\closure{\Omega})$, if $\abs{f(x)-f(y)} \leq \omega(\abs{x-y})$ for some modulus of continuity $\omega$, which is a continuous nondecreasing function satisfying 
\begin{equation*}
\omega(0) = 0, \qquad \int_0^1 \frac{\omega(t)}{t} \,dt < \infty.
\end{equation*}
The following result is from \cite{salothesis}.

\begin{thm} \label{thm:uniqueness}
Let $\Omega \subseteq \mR^n$ be a bounded open set with $C^1$ boundary, and $n \geq 3$. Suppose $A_1, A_2 \in C^{\mathrm{Dini}}(\closure{\Omega} ; \mC^n)$ and $V_1, V_2 \in L^{\infty}(\Omega ; \mC)$. Let $\Lambda_{A_1,V_1} = \Lambda_{A_2,V_2}$. Then $dA_1 = dA_2$ and $V_1 = V_2$.
\end{thm}

To prove this, one first uses boundary determination (Theorem \ref{thm:boundary}) and gauge invariance to get $A_1 = A_2$ on $\partial \Omega$, and then the vector fields are extended to a larger ball $B$ so that they coincide outside $\Omega$ . The DN maps on $\partial B$ coincide since the original DN maps do, and after another gauge transformation we may assume that $A_j \in C(\closure{B} ; \mC^n)$ and $D \cdot A_j = 0$ (this is where Dini continuity is required). The construction of CGO solutions in Theorem \ref{thm:cgo_pdo} can now be applied, and we may use the CGO solutions and arguments similar to those given in the reconstruction to show that $dA_1 = dA_2$ and $V_1 = V_2$.

We next consider stability results, which state that two sets of coefficients should be close if the DN maps are. It is known that many inverse problems are ill-posed, and for instance the Calder\'on problem only enjoys a logarithmic stability estimate. Logarithmic stability also holds for the magnetic Schr\"odinger equation. This was proved in \cite{tzou} along with a $\log\,\log$ stability estimate for partial data. We state the full data result in a slightly simplified form.

\begin{thm} \label{thm:stability}
Let $\Omega \subseteq \mR^n$ be a bounded open set with smooth boundary, and $n \geq 3$. Suppose $A_j \in W^{2,\infty}_c(\Omega ; \mR^n)$ and $V_j \in L^{\infty}(\Omega)$, and suppose $\norm{A_j}_{W^{2,\infty}} \leq M$, $\norm{V_j}_{L^{\infty}} \leq M$ for $j=1,2$. There exist $C > 0$ and $\varepsilon > 0$ such that 
\begin{multline*}
\norm{d(A_1-A_2)}_{H^{-1}} + \norm{V_1-V_2}_{H^{-1}} \leq C(\norm{\Lambda_{A_1,V_1}-\Lambda_{A_2,V_2}}^{1/2} \\
 + \abs{\log\,\norm{\Lambda_{A_1,V_1}-\Lambda_{A_2,V_2}}}^{-\varepsilon})
\end{multline*}
where $\norm{\,\cdot\,} = \norm{\,\cdot\,}_{H^{1/2}(\partial \Omega) \to H^{-1/2}(\partial \Omega)}$ on the right.
\end{thm}

A main step in the proof is to ensure that whenever the coefficients of $H_{A,V}$ satisfy a priori bounds, such as 
\begin{equation*}
\norm{A}_{W^{2,\infty}} \leq M, \quad \norm{V}_{L^{\infty}(\Omega)} \leq M, \quad A \text{ and } V \text{ supported in } B(0,M)
\end{equation*}
where $M > 0$, then the construction of CGO solutions in Theorem \ref{thm:cgo_pdo} goes through with constants only depending on $M$. This means that there are $C_M > 0$, $\varepsilon_M > 0$, so that the solutions exist for $\abs{\rho} \geq C_M$ and the remainder $\tilde{r}_{\rho}$ satisfies 
\begin{equation*}
\abs{\rho} \norm{\tilde{r}_{\rho}}_{L^2_{\delta}} + \norm{\tilde{r}_{\rho}}_{L^2_{\delta}} \leq C_M \abs{\rho}^{1-\varepsilon_M}.
\end{equation*}
Given this, one may insert suitable CGO solutions in the weak definition of $\Lambda_{A_1,V_1}-\Lambda_{A_2,V_2}$, and the Fourier transform of components of $d(A_1-A_2)$ will be controlled by quantities depending on $\norm{\Lambda_{A_1,V_1}-\Lambda_{A_2,V_2}}$ and $M$. For the electric potentials, one needs additional gauge transformations to control $\norm{A_1-A_2}_{L^{\infty}}$ by $\norm{d(A_1-A_2)}_{L^p}$ where $p > n$.

\section{CGO solutions and Carleman estimates} \label{sec:cgo_carleman}

Recently, a new method was introduced for constructing CGO solutions to the magnetic Schr\"odinger equation. The method is based on Carleman estimates, and it avoids the use of pseudodifferential operators. We will describe the main ideas of the construction here. In the next section, we discuss how these ideas can be used to prove partial data results for $H_{A,V}$.

If one is working in a bounded domain $\Omega \subseteq \mR^n$, the Sylvester-Uhlmann estimates in Proposition \ref{prop:sylvesteruhlmann} may be written as 
\begin{equation} \label{carleman_first}
h \norm{u} + h \norm{hDu} \leq C \norm{e^{\varphi/h} (-h^2 \Delta) e^{-\varphi/h} u}, \quad u \in C^{\infty}_c(\Omega),
\end{equation}
where $\norm{\,\cdot\,} = \norm{\,\cdot\,}_{L^2(\Omega)}$, $\varphi(x) = \alpha \cdot x$ with $\alpha \in \mR^n$ a unit vector, and $h > 0$ is small. Such weighted norm estimates with small parameter are called Carleman estimates. Now \eqref{carleman_first} is sufficient for constructing CGO solutions of the form \eqref{cgo_def_1} to $H_{0,V} u = 0$ in $\Omega$. This point of view was adopted by Bukhgeim and Uhlmann \cite{bukhgeimuhlmann}, who used \eqref{carleman_first} to give a partial data result for the operator $H_{0,V}$.

Kenig, Sj\"ostrand, and Uhlmann \cite{ksu} introduced new techniques in inverse boundary problems using the Carleman estimate approach. They consider more general limiting Carleman weights $\varphi$ for which estimates like \eqref{carleman_first} can be proved, both for $\varphi$ and $-\varphi$. The main examples are the linear weight $\varphi(x) = \alpha \cdot x$, and the logarithmic weight 
\begin{equation} \label{logweight_def}
\varphi(x) = \log\,\abs{x-x_0}, \quad x_0 \notin \closure{\Omega}.
\end{equation}
Approximate CGO solutions to $H_{0,V} u = 0$ are obtained from a WKB construction for the conjugated operator $e^{\varphi/h} h^2 H_{0,V} e^{-\varphi/h}$, and the approximate solutions are converted into exact CGO solutions by solving an inhomogeneous equation. The last step is made possible by the Carleman estimate, which also yields decay of the error term in the $L^2$ norm.

Proving the estimate \eqref{carleman_first} usually requires a convexity condition on the weight $\varphi$. Since the estimate also needs to hold for $-\varphi$, the limiting Carleman weights $\varphi$ will only satisfy a degenerate convexity condition. In \cite{ksu} a convexification argument was used to prove \eqref{carleman_first}: one replaces the limiting Carleman weight $\varphi$ by a convexified weight $\varphi_{\varepsilon}$, for instance 
\begin{equation} \label{convexified_weight}
\varphi_{\varepsilon} = \varphi + h \frac{\varphi^2}{2 \varepsilon}.
\end{equation}
One then proves the Carleman estimate for $\varphi_{\varepsilon}$, and uses this to deduce the estimate for the original weight $\varphi$.

The preceding discussion applied to CGO solutions of $H_{0,V} u = 0$. For the magnetic operator $H_{A,V}$ where $A$ is large, \eqref{carleman_first} does not immediately provide CGO solutions, since the first order term is not a small perturbation. However, if one uses the convexified weight $\varphi_{\varepsilon}$ then \eqref{carleman_first} can be proved also in the magnetic case, and this will imply the Carleman estimate for $H_{A,V}$ with the original weight $\varphi$. In other words, convexification makes the procedure more robust, so that one can handle first order perturbations. This method was used in \cite{dksu} to construct CGO solutions to $H_{A,V} u = 0$, and to give a partial data result for the magnetic Schr\"odinger equation.

The results in \cite{dksu} assume that the magnetic potential is $C^2$. In \cite{knudsensalo}, the results were extended to H\"older continuous potentials. The problem in passing from $C^2$ to $C^{\varepsilon}$ magnetic potentials is that the regularity of CGO solutions deteriorates. In particular, the solutions will not be in $H^1(\Omega)$. However, by combining Carleman estimates with the pseudodifferential conjugation argument, we obtain CGO solutions in $H^1(\Omega)$ even for H\"older coefficients. One of the main points of \cite{knudsensalo} is an extension of the pseudodifferential conjugation argument to logarithmic Carleman weights.

In the rest of this section, we give some more details on the construction of CGO solutions using Carleman estimates. The main result will be the following. This was proved for $A \in C^2$ in \cite{dksu}, and for $A \in C^{\varepsilon}$ in \cite{knudsensalo}.

\begin{thm} \label{thm:cgo_carleman}
Let $A \in C^{\varepsilon}(\closure{\Omega} ; \mC^n)$ and $\tilde{V} \in L^n(\Omega)$. Let $\varphi$ and $\psi$ be defined by \eqref{logweight_def} and \eqref{psidef}, respectively. Then for $h$ small there is an $H^1(\Omega)$ solution $u = e^{-\frac{1}{h}(\varphi+i\psi)}(a+r)$ of the equation $H_{A,V} u = 0$ in $\Omega$, where $a$ is given by \eqref{wkb_amplitude_def}. One has the norm estimates
\begin{eqnarray*}
 & \norm{\partial^{\alpha} a}_{L^{\infty}(\Omega)} = O(h^{-\sigma\abs{\alpha}}), & \\
 & \norm{r}_{L^2(\Omega)} + \norm{h\nabla r}_{L^2(\Omega)} = O(h^{\sigma\varepsilon}) & 
\end{eqnarray*}
where $\sigma > 0$ is small.
\end{thm}

We follow \cite{ksu} and \cite{dksu} and assume first that $A \in C^2$. Let $\Omega, \tilde{\Omega} \subseteq \mR^n$ be bounded domains with $\closure{\Omega} \subseteq \tilde{\Omega}$, and let $\varphi \in C^{\infty}(\tilde{\Omega} ; \mR)$ with $\nabla \varphi \neq 0$ in $\tilde{\Omega}$. We want to find conditions on $\varphi$ so that \eqref{carleman_first} holds for $\varphi$ and $-\varphi$. Consider the conjugated operator 
\begin{equation*}
P_{\varphi} = e^{\varphi/h}(-h^2 \Delta)e^{-\varphi/h}.
\end{equation*}
We would like to show $\norm{P_{\varphi} u} \geq c h\norm{u}$ when $u \in C_c^{\infty}(\Omega)$. To do this, we write $P_{\varphi} = A + iB$ where $A, B$ are self-adjoint, and note that 
\begin{equation*}
\norm{P_{\varphi}u}^2 = (P_{\varphi}u|P_{\varphi}u) = \norm{Au}^2 + \norm{Bu}^2 + (i[A,B]u|u),
\end{equation*}
where $[A,B] = AB-BA$ is the commutator. To get a bound from below, one would like a positivity condition for $i[A,B]$. The principal symbol (in semiclassical Weyl quantization, see Section \ref{sec:cgo_pdo}) of $i[A,B]$ is $h\{a,b\}$, where $\{a,b\} = \nabla_{\xi} a \cdot \nabla_x b - \nabla_x a \cdot \nabla_{\xi} b$ is the Poisson bracket, and 
\begin{equation*}
a(x,\xi) = \xi^2 - (\nabla \varphi)^2, \quad b(x,\xi) = 2\nabla \varphi \cdot \xi \quad (x \in \tilde{\Omega}, \xi \in \mR^n).
\end{equation*}
The positivity condition for $i[A,B]$ at symbol level is 
\begin{equation*}
\{a,b\} \geq 0 \quad \text{when } a = b = 0.
\end{equation*}
Since $[A,B]$ changes sign when $\varphi$ is replaced by $-\varphi$, one also needs $\{a,b\} \leq 0$ when $a = b = 0$. This leads to the following condition.

\begin{definition}
Let $\varphi$ be as above. We say that $\varphi$ is a limiting Carleman weight (for the Laplacian) if 
\begin{equation} \label{limiting_carleman}
\{a,b\} = 0 \quad \text{when } a = b = 0.
\end{equation}
\end{definition}

The condition \eqref{limiting_carleman} is used in many ways in the construction of CGO solutions and determination of coefficients. However, the condition is also quite restrictive.

Given \eqref{limiting_carleman}, we want to prove a Carleman estimate. We first convexify the weight and replace $\varphi$ by $\varphi_{\varepsilon}$, defined as in \eqref{convexified_weight}. If $A_{\varepsilon}$ and $B_{\varepsilon}$ are the operators corresponding to $A$ and $B$, one still has 
\begin{equation} \label{epsilon_carlemanproof}
\norm{P_{\varphi_{\varepsilon}}u}^2 = \norm{A_{\varepsilon}u}^2 + \norm{B_{\varepsilon}u}^2 + (i[A_{\varepsilon},B_{\varepsilon}]u|u).
\end{equation}
Now, if $\varepsilon$ is small enough, one can prove the lower bound 
\begin{equation*}
(i[A_{\varepsilon},B_{\varepsilon}]u|u) \geq \frac{c h^2}{\varepsilon} \norm{u}^2 - \frac{1}{2} \norm{A_{\varepsilon}u}^2 - \frac{1}{2} \norm{B_{\varepsilon}u}^2.
\end{equation*}
Note that the first term is positive. This yields the estimate 
\begin{equation*}
\frac{h}{\sqrt{\varepsilon}} (\norm{u} + \norm{hDu}) \leq C \norm{e^{\varphi_{\varepsilon}/h} (-h^2\Delta) e^{-\varphi_{\varepsilon}/h} u}.
\end{equation*}
By choosing $\varepsilon$ small enough, we can replace $-h^2 \Delta$ by $h^2 H_{A,V}$. Finally, since $e^{\varphi_{\varepsilon}/h} = e^{\varphi/h} m$ where $m$ is bounded from above and below in $\closure{\Omega}$, we arrive at the Carleman estimate 
\begin{equation} \label{carleman_hav}
h (\norm{u} + \norm{hDu}) \leq C \norm{e^{\varphi/h} h^2 H_{A,V} e^{-\varphi/h} u}, \quad u \in C_c^{\infty}(\Omega).
\end{equation}
Shifting the estimate to a different Sobolev index and using the Hahn-Banach theorem, one may show that this a priori estimate implies solvability for an inhomogeneous equation.

\begin{prop} \label{prop:cgocarleman_inhomog}
\cite{dksu} Let $\varphi$ be a limiting Carleman weight in $\tilde{\Omega}$, let $A \in C^1(\closure{\Omega})$, and let $V \in L^{\infty}(\Omega)$. If $h$ is small enough, then for any $f \in L^2(\Omega)$, the equation 
\begin{equation*}
e^{\varphi/h} h^2 H_{A,V} e^{-\varphi/h} r = f \quad \text{in } \Omega
\end{equation*}
has a solution $r \in H^1(\Omega)$, which satisfies $h\norm{r} + h\norm{hDr} \leq C \norm{f}$.
\end{prop}

Proposition \ref{prop:cgocarleman_inhomog} will be used to convert an approximate CGO solution to an exact one. Given a limiting Carleman weight $\varphi$, the CGO solution $u = e^{-\varphi/h} v$ will satisfy 
\begin{equation} \label{carleman_approx_eq}
e^{\varphi/h} h^2 H_{A,V} e^{-\varphi/h} v = 0.
\end{equation}
We find $v$ by a WKB construction in the form 
\begin{equation*}
v = e^{-\frac{1}{h} i\psi} (a+r)
\end{equation*}
where $\psi \in C^{\infty}(\closure{\Omega} ; \mR)$ is a phase function, $a$ is an amplitude, and $r$ is a correction term. Inserting the WKB ansatz in \eqref{carleman_approx_eq}, and grouping the terms in powers of $h$, one obtains 
\begin{equation*}
(-(\nabla \rho)^2 + ih[\nabla \rho \circ D + D \circ \nabla \rho + 2\nabla\rho \cdot A] + h^2 H_{A,V})(a+r) = 0,
\end{equation*}
where $\rho = \varphi+i\psi$. This will be satisfied if $\psi$, $a$, and $r$ are chosen so that 
\begin{eqnarray}
 & (\nabla \rho)^2 = 0, & \label{eikonal} \\
 & (\nabla \rho \circ D + D \circ \nabla \rho + 2\nabla\rho \cdot A)a = 0, & \label{transport} \\
 & e^{\varphi/h} h^2 H_{A,V} e^{-\varphi/h}(e^{-i\frac{\psi}{h}} r) = -h^2 e^{-i\frac{\psi}{h}} H_{A,V} a. & \label{correction}
\end{eqnarray}
The first equation may be written as 
\begin{equation*}
(\nabla \psi)^2 = (\nabla \varphi)^2, \quad \nabla \varphi \cdot \nabla \psi = 0.
\end{equation*}
This is an eikonal equation for $\psi$. If $\psi$ is known, \eqref{transport} is a transport equation for $a$, and when $\psi$ and $a$ have been found, the last equation for $r$ can be solved by Proposition \ref{prop:cgocarleman_inhomog}.

We now choose $\varphi$ as in \eqref{logweight_def}, and indicate how to solve \eqref{eikonal} -- \eqref{transport} following \cite{ksu}, \cite{dksu}. It is known that distance functions are solutions of the eikonal equation, and one may take 
\begin{equation} \label{psidef}
\psi(x) = \mathrm{dist}_{S^{n-1}}\Big( \frac{x-x_0}{\abs{x-x_0}}, \omega \Big),
\end{equation}
where $\omega \in S^{n-1}$ is chosen so that $\psi$ is smooth near $\closure{\Omega}$. If one chooses coordinates so that $x_0 = 0$, $\closure{\Omega}$ lies in $\{x_n > 0\}$, $\omega = e_1$, and $x = (x_1,r\theta)$ with $r > 0$ and $\theta \in S^{n-1}$, then 
\begin{equation*}
\varphi + i\psi = \log\,z
\end{equation*}
where $z = x_1 + ir$ is a complex variable. In these new coordinates, the transport equation \eqref{transport} becomes 
\begin{equation*}
(\partial_{\bar{z}} + \frac{i}{2}(e_1+i e_r) \cdot A(z,\theta) - \frac{n-2}{2(z-\bar{z})})a(z,\theta) = 0
\end{equation*}
where $e_r = (0,\theta)$. This has a solution 
\begin{equation} \label{wkb_amplitude_def}
a(z,\theta) = (z-\bar{z})^{\frac{2-n}{2}} e^{i\Phi}
\end{equation}
where $\Phi(z,\theta) = -\frac{1}{2} \partial_{\bar{z}}^{-1}((e_1+ie_r) \cdot A(z,\theta))$ is obtained as a Cauchy transform. This ends the proof of Theorem \ref{thm:cgo_carleman} in the case $A \in C^2$.

If $A$ is only $C^{\varepsilon}$, the main step is to prove Proposition \ref{prop:cgocarleman_inhomog} under this assumption. We follow \cite{knudsensalo} and write, as in \eqref{conjugated_identity_pdo},
\begin{equation*}
e^{\varphi/h} h^2 H_{A,V} e^{-\varphi/h} = P + hQ + h^2 \tilde{V},
\end{equation*}
where $P = e^{\varphi/h} (-h^2 \Delta) e^{-\varphi/h}$ and $Q = e^{\varphi/h} (2A \cdot hD) e^{-\varphi/h}$. The free equation $Pr = f$ can be solved by Proposition \ref{prop:cgocarleman_inhomog}, and one wishes to use a perturbation argument to solve the full equation. Again, the first order term $Q$ is not a small perturbation. We split $Q = Q^{\sharp} + Q^{\flat}$ where $Q^{\sharp} = e^{\varphi/h} (2A^{\sharp} \cdot hD) e^{-\varphi/h}$, with $A^{\sharp}$ a $h$-dependent smooth approximation as in \eqref{asharp_def}, with $\delta = h^{\sigma}$ and $\sigma > 0$ small. Then $Q^{\flat}$ will be a small perturbation which can be inverted by Neumann series, and the smooth part $Q^{\sharp}$ can be conjugated away using pseudodifferential operators as in Proposition \ref{prop:conjugation}.

\begin{prop} \label{prop:psdoconjugation_carleman}
There exist $c, \tilde{c}, r \in S^0_{\sigma}$ so that 
\begin{equation*}
(P+hQ^{\sharp})C = \tilde{C}P + h^{2-2\sigma} R \quad \text{in } \Omega.
\end{equation*}
Both $C$ and $\tilde{C}$ are elliptic, in the sense that $c$ and $\tilde{c}$ are nonvanishing for small $h$.
\end{prop}

The proof is similar to that of Proposition \ref{prop:conjugation}. Here one does not need global solutions, so it is enough to use cutoffs to restrict to bounded domains. However, the transport equation \eqref{hp_dbar_eq} does not have constant coefficients anymore, and more work is needed to show that it can be solved in a full neighborhood of $p^{-1}(0)$. To do this, we use ideas from \cite{duistermaathormander} and find a smooth function $m$ so that instead of \eqref{limiting_carleman}, one has $\{ma,mb\} = 0$ in a neighborhood of $a = b = 0$. The flows of the Hamilton vector fields $H_{ma}$ and $H_{mb}$ will then commute, and one can compute a change of coordinates in a neighborhood of $p^{-1}(0)$ which reduces $H_{mp}$ to $\partial_{y_1} + i\partial_{y_2}$. After this, the equation \eqref{hp_dbar_eq} is easy to solve. The details are given in \cite{knudsensalo}.

\section{Partial data results} \label{sec:partial}

Let $\Omega \subseteq \mR^n$, $n \geq 3$, be a bounded simply connected $C^{\infty}$ domain with connected boundary. We want to discuss the unique determination of the magnetic field $dA$ and electric potential $V$, if one knows the boundary measurements $\Lambda_{A,V}$ on a subset of the boundary $\partial \Omega$.

Let $x_0 \in \mR^n \smallsetminus \closure{\mathrm{ch}(\Omega)}$, where $\mathrm{ch}(\Omega)$ is the convex hull of $\Omega$. We define the front face of $\partial \Omega$ relative to $x_0$ by 
\begin{equation*}
F(x_0) = \{ x \in \partial \Omega \colon (x-x_0) \cdot \nu(x) \leq 0 \},
\end{equation*}
and we take $\tilde{F}$ to be an open neighborhood of $F(x_0)$ in $\partial \Omega$. Also, let $A_{\mathrm{tan}}$ be the tangential component defined in \eqref{tangential_component}.

The following result was proved for $C^2$ magnetic potentials in \cite{dksu}, and it was extended to H\"older continuous potentials in \cite{knudsensalo}.

\begin{thm} \label{thm:partialdata}
Let $A_j \in C^{\varepsilon}(\closure{\Omega} ; \mC^n),\; \varepsilon > 0,$ and $V_j \in L^{\infty}(\Omega)$ for $j = 1,2$. Also assume that $0$ is not a Dirichlet eigenvalue of $H_{A_j,V_j}$ in $\Omega$. If
\begin{equation*}
\Lambda_{A_1,V_1} f|_{\tilde{F}} = \Lambda_{A_2,V_2} f|_{\tilde{F}} \quad \text{for all } f \in H^{1/2}(\partial \Omega),
\end{equation*}
then $dA_1 = dA_2$ in $\Omega$ and $V_1 = V_2$ in $\Omega$. Also, $(A_1)_{\mathrm{tan}} = (A_2)_{\mathrm{tan}}$ on $\partial \Omega$.
\end{thm}

We describe the proof in the case where $A_j \in C^2(\closure{\Omega} ; \mC^n)$ with $A_j \cdot \nu = 0$ on $\partial \Omega$, following \cite{dksu}. See \cite{knudsensalo} for the $C^{\varepsilon}$ case and some further details. We start by choosing $u_1$ and $u_2$ to be CGO solutions, given by Theorem \ref{thm:cgo_carleman}, to the equations $H_{A_1,V_1} u_1 = 0$ and $H_{\bar{A}_2,\bar{V}_2} u_2 = 0$. They will be of the form 
\begin{align*}
u_1 &= e^{\frac{1}{h}(\varphi+i\psi)}(a_1+r_1), \\
u_2 &= e^{\frac{1}{h}(-\varphi+i\psi)}(a_2+r_2),
\end{align*}
where $\varphi$ is the logarithmic weight \eqref{logweight_def}, $\psi$ is the distance function \eqref{psidef}, and $a_j$ are amplitudes solving transport equations. One reason for working with limiting Carleman weights is that $u_1$ and $u_2$ can be chosen with opposite signs in front of $\varphi$. With this choice, the terms $e^{\pm \frac{1}{h} \varphi}$ which grow as $h \to 0$ cancel in the product $u_1 \bar{u}_2$.

Take $\tilde{u}_2$ to be the solution to $H_{A_2,V_2} \tilde{u}_2 = 0$ with $\tilde{u}_2 = u_1$ on $\partial \Omega$, and let $u = u_1 - \tilde{u}_2$. Integration by parts shows 
\begin{equation} \label{partial_green}
(H_{A_2,V_2} u | u_2) = -(\partial_{\nu} u | u_2)_{\partial \Omega}.
\end{equation}
Since the DN maps coincide on $\tilde{F}$, and since $\partial_{\nu} u = (\Lambda_{A_1,V_1} - \Lambda_{A_2,V_2}) u_1$, one sees that $\partial_{\nu} u = 0$ on $\tilde{F}$. Thus the right hand side of \eqref{partial_green} reduces to an integral over $\partial \Omega \smallsetminus \tilde{F}$, and one would like to estimate $\partial_{\nu} u$ in this set.

As shown in \cite{bukhgeimuhlmann}, one can estimate $\partial_{\nu} u$ by using a Carleman estimate with boundary terms. This is an extension of \eqref{carleman_hav} to functions $v$ which vanish on $\partial \Omega$ but are not necessarily compactly supported in $\Omega$. The estimate has the form 
\begin{multline*}
\sqrt{h} \norm{\sqrt{\partial_{\nu} \varphi} e^{-\frac{\varphi}{h}} \partial_{\nu} v}_{L^2(\partial \Omega_{+})} + \norm{e^{-\frac{\varphi}{h}} v} + \norm{e^{-\frac{\varphi}{h}} h \nabla v} \\
 \leq Ch \norm{e^{-\frac{\varphi}{h}} H_{A_2,V_2} v} + C\sqrt{h} \norm{\sqrt{-\partial_{\nu} \varphi} e^{-\frac{\varphi}{h}} \partial_{\nu} v}_{L^2(\partial \Omega_{-})}.
\end{multline*}
Here $\partial \Omega_{\pm} = \{ x \in \partial \Omega \colon \pm \partial_{\nu} \varphi(x) \geq 0\}$ are the front and back sides of $\partial \Omega$. This estimate can be applied to $u$, and one gets a bound for $\norm{e^{-\frac{\varphi}{h}} \partial_{\nu} u}_{L^2(\partial \Omega \smallsetminus \tilde{F})}$. Consequently 
\begin{equation*}
h (\partial_{\nu} u | u_2)_{\partial \Omega} \to 0
\end{equation*}
as $h \to 0$. Computing the limit of $h$ times the left hand side of \eqref{partial_green}, and using the explicit forms for $u_1$ and $u_2$, yields 
\begin{equation} \label{partial_firstorthogonal}
\int_{\Omega} \nabla (\varphi+i\psi) \cdot (A_1-A_2) a \,dx = 0,
\end{equation}
where $a = \lim_{h \to 0} a_1 \bar{a}_2$. The expression on the left may be thought of as a counterpart to \eqref{htav_limit}, and we will see below that \eqref{partial_firstorthogonal} asserts the vanishing of a certain nonlinear Radon transform of $d(A_1-A_2)$. This will imply $dA_1 = dA_2$.

We switch to the complex notation used in Section \ref{sec:cgo_carleman}. Thus, $z = x_1 +ir$ where $x = (x_1,r\theta)$. We also write $P_{\theta}$ for the two-plane consisting of points $(x_1,r\theta)$ with $\theta$ fixed, and $\Omega_{\theta} = \Omega \cap P_{\theta}$. The formula \eqref{partial_firstorthogonal} becomes 
\begin{equation} \label{partial_secondorthogonal}
\int_{S^{n-2}} \Big( \int_{\Omega_{\theta}} \frac{1}{z} (e_1+ie_r) \cdot (A_1-A_2)(z,\theta) e^{i\Phi} \,d\bar{z} \wedge dz \Big) \,d\theta = 0.
\end{equation}
Here $\Phi$ arises from the amplitudes $a_1$ and $a_2$, and will satisfy 
\begin{equation} \label{phitransport}
\partial_{\bar{z}} \Phi(z,\theta) = -\frac{1}{2}(e_1+i e_r) \cdot (A_1-A_2)(z,\theta).
\end{equation}
We need a slightly more general version of \eqref{partial_secondorthogonal}. From the transport equation for $a_1$, one sees that $a_1$ can be replaced by $a_1 g_1$ where $\partial_{\bar{z}} g_1(z,\theta) = 0$. Then \eqref{partial_secondorthogonal} holds with $e^{i\Phi}$ replaced with $e^{i\Phi} g_1$. Choosing $g_1(z,\theta) = z g(z) \tilde{g}(\theta)$ with $\partial_{\bar{z}} g = 0$, and varying $\tilde{g}$, we see that for almost every $\theta$ one has 
\begin{equation} \label{secondtransform}
\int_{\Omega_{\theta}} (e_1+ie_r) \cdot (A_1-A_2)(z,\theta) e^{i\Phi(z,\theta)} g(z) \,d\bar{z} \wedge dz = 0.
\end{equation}

We will show that \eqref{secondtransform} remains true with $e^{i\Phi}$ replaced with $1$, which corresponds to converting a nonlinear Radon transform to the usual one. This will be done by ''taking a holomorphic logarithm''. The procedure has similarities with \cite{sun}. The first step is to use the equation \eqref{phitransport} and integrate by parts in \eqref{secondtransform}, which yields 
\begin{equation} \label{eiphi_orthogonal}
\int_{\partial \Omega_{\theta}} e^{i\Phi} g(z) \,dz = 0
\end{equation}
for $g$ holomorphic in $\Omega_{\theta}$. This orthogonality condition means that $e^{i\Phi}|_{\partial \Omega_{\theta}}$ is the boundary value of some holomorphic function $F \in C(\closure{\Omega}_{\theta})$. It can be shown that $F$ is nonvanishing and has a holomorphic logarithm $G \in C(\closure{\Omega}_{\theta})$. Taking $g = \frac{G}{e^G}$ in \eqref{eiphi_orthogonal}, and noting that $G = i\Phi$ on $\partial \Omega$ up to constant, one gets 
\begin{equation*}
\int_{\partial \Omega_{\theta}} \Phi \,dz = 0.
\end{equation*}
Integrating by parts and using \eqref{phitransport} again gives 
\begin{equation*}
\int_{\Omega_{\theta}} (e_1+i e_r) \cdot (A_1-A_2)(z,\theta) \,d\bar{z} \wedge dz = 0.
\end{equation*}
Going back to $x$-coordinates, and varying $x_0$ and $\omega$ slightly, this will imply 
\begin{equation} \label{partial_radonvanishing}
\int_{P \cap \Omega} \xi \cdot (A_1-A_2) \,dS = 0
\end{equation}
for all two-planes $P$ with $d((0,e_1),T(P)) < \delta$, and all $\xi \in P$.

Finally, we indicate how \eqref{partial_radonvanishing} implies $dA_1 = dA_2$ in the case $n = 3$. The argument uses analytic microlocal analysis. By approximation, one may assume $A_j \in C^{\infty}_c(\mR^3)$ in \eqref{partial_radonvanishing}. Fixing $\xi, \eta \in \mR^3$ and writing $A = A_1-A_2$, we obtain from \eqref{partial_radonvanishing} 
\begin{equation*}
\int_P \br{dA,\xi \wedge \eta} \,dS = 0
\end{equation*}
for all $P$ with $d((0,e_1),T(P)) < \delta$. By the microlocal Helgason support theorem, the vanishing of the Radon transform at the planes $P$ implies the microlocal analyticity of $\br{dA,\xi \wedge \eta}$ on the conormal bundle of every such $P$. Consequently, if $\br{dA,\xi \wedge \eta}$ is nonzero, then it is microlocally analytic at some point of the normal set of $\supp(\br{dA,\xi \wedge \eta})$. But this contradicts the microlocal Holmgren theorem, so one must have $dA_1 = dA_2$. We refer to \cite{dksu} for more details, and for the case of electric potentials.


\begin{thebibliography}{10}

\bibitem{brownboundary}
R.~M. Brown, \emph{{Recovering conductivity at the boundary from the Dirichlet
  to Neumann map: a pointwise result}}, J. Inverse Ill-Posed Probl. \textbf{9}
  (2001), 567--574.

\bibitem{brownsalo}
R.~M. Brown and M.~Salo, \emph{Identifiability at the boundary for first-order
  terms}, Appl. Anal. \textbf{85} (2006), no.~6-7, 735--749.

\bibitem{bukhgeimuhlmann}
A.~L. Bukhgeim and G.~Uhlmann, \emph{{Recovering a potential from partial
  Cauchy data}}, Comm. PDE \textbf{27} (2002), 653--668.

\bibitem{calderon}
A.~P. Calder{\'o}n, \emph{On an inverse boundary value problem}, Seminar on
  Numerical Analysis and its Applications to Continuum Physics, Soc. Brasileira
  de Matem{\'a}tica, R{\'i}o de Janeiro, 1980.

\bibitem{dimassisjostrand}
M.~Dimassi and J.~Sj{\"o}strand, \emph{Spectral asymptotics in the
  semi-classical limit}, London Mathematical Society Lecture Note Series 268,
  Cambridge University Press, 1999.

\bibitem{dksu}
D.~Dos Santos~Ferreira, C.~E. Kenig, J.~Sj{\"o}strand, and G.~Uhlmann,
  \emph{{Determining a magnetic Schr{\"o}dinger operator from partial Cauchy
  data}}, Comm. Math. Phys. (to appear), \mbox{arXiv:math.AP/0601466}.

\bibitem{duistermaathormander}
J.~J. Duistermaat and L.~H{\"o}rmander, \emph{{Fourier integral operators.
  II}}, Acta Math. \textbf{128} (1972), 183--269.

\bibitem{eskin}
G.~Eskin, \emph{{Global uniqueness in the inverse scattering problem for the
  Schr{\"o}dinger operator with external Yang-Mills potentials}}, Comm. Math.
  Phys. \textbf{222} (2001), no.~3, 503--531.

\bibitem{eskinralston}
G.~Eskin and J.~Ralston, \emph{{Inverse scattering problem for the
  Schr{\"o}dinger equation with magnetic potential at a fixed energy}}, Comm.
  Math. Phys. \textbf{173} (1995), 199--224.

\bibitem{eskinralstonproc}
\bysame, \emph{{Inverse scattering problems for Schr{\"o}dinger operators with
  magnetic and electric potentials}}, Inverse problems in wave propagation
  ({Guy Chavent et al.}, ed.), IMA Volumes in Mathematics and Applications,
  no.~90, Springer, 1997, pp.~147--166.

\bibitem{eskinralston_elasticity}
\bysame, \emph{On the inverse boundary value problem for linear isotropic
  elasticity}, Inverse Problems \textbf{18} (2002), 907--921.

\bibitem{heckliwang}
H.~Heck, X.~Li, and J.-N. Wang, \emph{Identification of viscosity in an
  incompressible fluid}, Indiana Univ. Math. J. (to appear).

\bibitem{jkexp}
D.~Jerison and C.~Kenig, \emph{{Boundary value problems on Lipschitz domains}},
  Studies in partial differential equations (Walter Littman, ed.), MAA Studies
  in Mathematics, no.~23, Mathematical Association of America, 1982, pp.~1--68.

\bibitem{kenigponcevega}
C.~E. Kenig, G.~Ponce, and L.~Vega, \emph{{Smoothing effects and local
  existence theory for the generalized nonlinear Schr\"odinger equations}},
  Invent. Math. \textbf{134} (1998), 489--545.

\bibitem{ksu}
C.~E. Kenig, J.~Sj{\"o}strand, and G.~Uhlmann, \emph{{The Calder\'on problem
  with partial data}}, Ann. of Math. (to appear), \mbox{arXiv:math.AP/0405486}.

\bibitem{knudsensalo}
K.~Knudsen and M.~Salo, \emph{Determining nonsmooth first order terms from
  partial boundary measurements}, Inverse Problems and Imaging (to appear),
  \mbox{arXiv:math.AP/0609133}.

\bibitem{kohnvogelius1}
R.~Kohn and M.~Vogelius, \emph{Determining conductivity by boundary
  measurements}, Comm. Pure Appl. Math. \textbf{37} (1984), 289--298.

\bibitem{lidirac}
X.~Li, \emph{{On the inverse problem for the Dirac operator}}, preprint.

\bibitem{mcdowall}
S.~R. McDowall, \emph{An electromagnetic inverse problem in chiral media},
  Trans. Amer. Math. Soc. \textbf{352} (2000), no.~7, 2993--3013.

\bibitem{nachman}
A.~Nachman, \emph{Reconstructions from boundary measurements}, Ann. of Math.
  \textbf{128} (1988), 531--587.

\bibitem{nakamurasunuhlmann}
G.~Nakamura, Z.~Sun, and G.~Uhlmann, \emph{{Global identifiability for an
  inverse problem for the Schr{\"o}dinger equation in a magnetic field}}, Math.
  Ann. \textbf{303} (1995), 377--388.

\bibitem{nakamuratsuchida}
G.~Nakamura and T.~Tsuchida, \emph{{Uniqueness for an inverse boundary value
  problem for Dirac operators}}, Comm. PDE \textbf{25} (2000), 1327--1369.

\bibitem{nakamurauhlmann}
G.~Nakamura and G.~Uhlmann, \emph{Global uniqueness for an inverse boundary
  problem arising in elasticity}, Invent. Math. \textbf{118} (1994), 457--474.

\bibitem{nakamurauhlmannerratum}
\bysame, \emph{Erratum: Global uniqueness for an inverse boundary value problem
  arising in elasticity}, Invent. Math \textbf{152} (2003), 205--207.

\bibitem{novikovkhenkin}
R.G. Novikov and G.M. Khenkin, \emph{{The $\bar{\partial}$-equation in the
  multidimensional inverse scattering problem}}, Russ. Math. Surv. \textbf{42}
  (1987), 109--180.

\bibitem{panchenko}
A.~Panchenko, \emph{{An inverse problem for the magnetic Schr{\"o}dinger
  equation and quasi-exponential solutions of nonsmooth partial differential
  equations}}, Inverse Problems \textbf{18} (2002), 1421--1434.

\bibitem{salothesis}
M.~Salo, \emph{{Inverse problems for nonsmooth first order perturbations of the
  Laplacian}}, Ann. Acad. Sci. Fenn. Math. Diss. \textbf{139} (2004), 67 pp.

\bibitem{saloreconstruction}
\bysame, \emph{{Semiclassical pseudodifferential calculus and the
  reconstruction of a magnetic field}}, Comm. PDE (to appear),
  \mbox{arXiv:math.AP/0602290}.

\bibitem{sun}
Z.~Sun, \emph{{An inverse boundary value problem for Schr{\"o}dinger operators
  with vector potentials}}, Trans. Amer. Math. Soc. \textbf{338} (1993), no.~2,
  953--969.

\bibitem{sylvesteruhlmann}
J.~Sylvester and G.~Uhlmann, \emph{A global uniqueness theorem for an inverse
  boundary value problem}, Ann. of Math. \textbf{125} (1987), 153--169.

\bibitem{sylvesteruhlmannboundary}
\bysame, \emph{Inverse boundary value problems at the boundary - continuous
  dependence}, Comm. Pure Appl. Math. \textbf{41} (1988), no.~2, 197--219.

\bibitem{tolmasky}
C.~Tolmasky, \emph{{Exponentially growing solutions for nonsmooth first-order
  perturbations of the Laplacian}}, SIAM J. Math. Anal. \textbf{29} (1998),
  116--133.

\bibitem{tzou}
L.~Tzou, \emph{{Stability estimates for coefficients of magnetic
  Schr{\"o}dinger equation from full and partial boundary measurements}},
  \mbox{arXiv:math.AP/0602147}.

\bibitem{uhlmanndevelopments}
G.~Uhlmann, \emph{{Developments in inverse problems since Calder{\'o}n's
  foundational paper}}, {Harmonic analysis and partial differential equations.
  Essays in honor of Alberto P. Calder{\'o}n}, Univ. Chicago Press, Chicago,
  1999, pp.~295--345.

\bibitem{uhlmannselecta}
\bysame, \emph{{Commentary on Alberto P. Calder{\'o}n's paper: On an inverse
  boundary value problem}}, Selecta (to appear) (A.~Bellow, C.~E. Kenig, and
  P.~Malliavin, eds.).

\end{thebibliography}

\providecommand{\bysame}{\leavevmode\hbox to3em{\hrulefill}\thinspace}
\providecommand{\href}[2]{#2}

\end{document}